\documentclass[%
  onecolumn
   , hidempi
]{mpi2015-cscpreprint}
\makeatletter
\DeclareOldFontCommand{\rm}{\normalfont\rmfamily}{\mathrm}
\DeclareOldFontCommand{\sf}{\normalfont\sffamily}{\mathsf}
\DeclareOldFontCommand{\tt}{\normalfont\ttfamily}{\mathtt}
\DeclareOldFontCommand{\bf}{\normalfont\bfseries}{\mathbf}
\DeclareOldFontCommand{\it}{\normalfont\itshape}{\mathit}
\DeclareOldFontCommand{\sl}{\normalfont\slshape}{\@nomath\sl}
\DeclareOldFontCommand{\sc}{\normalfont\scshape}{\@nomath\sc}
\makeatother

\usepackage[american]{babel}

\usepackage{graphicx}
\usepackage{amssymb}
\usepackage{amsthm}
\usepackage[mathscr]{eucal}
\usepackage{enumitem}

\numberwithin{equation}{section}
\numberwithin{figure}{section}
\numberwithin{table}{section}

\usepackage[software,hardware]{mymacros}
\usepackage{import}
\usepackage{xifthen}
\usepackage{pdfpages}
\usepackage{transparent}
\usepackage{cleveref}
\usepackage{sidecap}    
\usepackage{floatrow}
\usepackage{todonotes}
\usepackage{algorithm}
\usepackage[noend]{algpseudocode}

\makeatletter
\def\algbackskip{\hskip-\ALG@thistlm}
\makeatother

\usepackage{caption}
\usepackage{subcaption}
\usepackage{arydshln}

\definecolor{lightgray}{gray}{0.9}
\usepackage{color}
\definecolor{bluegreen}{rgb}{0.0, 0.87, 0.87}
\usetikzlibrary{fadings,shapes.arrows,shadows}   
\usetikzlibrary{shapes.multipart}

\tikzfading[name=arrowfading, top color=transparent!0, bottom color=transparent!95]
\tikzset{arrowfill/.style={top color= black!20, bottom color=green!50!black, general shadow={fill=black, shadow yshift=-0.8ex, path fading=arrowfading}}}
\tikzset{arrowstyle/.style={draw=black,arrowfill, single arrow,minimum height=#1, single arrow,
		single arrow head extend=.2cm,}}

\newtheorem{remark}{Remark}

\newtheorem{theorem}{Theorem}[section]

\newtheorem{lemma}[theorem]{Lemma}
\newtheorem{problem}[theorem]{Problem}

\newcommand{\quadembs}{\texttt{quad-embeds}}
\newcommand{\linearembs}{\texttt{linear-embeds}}
\newcommand{\quadopInf}{\texttt{quad-OpInf}}

\begin{document}
  

\title{Generalized Quadratic-Embeddings for Nonlinear Dynamics using Deep Learning}
  
\author[$\ast$]{Pawan Goyal}
\affil[$\ast$]{Max Planck Institute for Dynamics of Complex Technical Systems, 39106 Magdeburg, Germany.\authorcr
  \email{goyalp@mpi-magdeburg.mpg.de}, \orcid{0000-0003-3072-7780}
}
  
\author[$\dagger\ddagger$]{Peter Benner}
\affil[$\dagger$]{Max Planck Institute for Dynamics of Complex Technical Systems, 39106 Magdeburg, Germany.\authorcr
  \email{benner@mpi-magdeburg.mpg.de}, \orcid{0000-0003-3362-4103}
}
\affil[$\ddagger$]{Otto von Guericke University,  Universit\"atsplatz 2, 39106 Magdeburg, Germany\authorcr
  \email{peter.benner@ovgu.de}
}
  
\shorttitle{Generalized Quadratic Embeddings for Nonlinear Dynamics}
\shortauthor{P. Goyal, P. Benner}
\shortdate{}
  

  
\abstract{%
The engineering design process often relies on mathematical modeling that can describe the underlying dynamic behavior.  
In this work, we present a data-driven methodology for modeling the dynamics of nonlinear systems. To simplify this task, we aim to identify a coordinate transformation that allows us to represent the dynamics of nonlinear systems using a common, simple model structure. The advantage of a common simple model is that customized design tools developed for it can be applied to study a large variety of nonlinear systems. 
The simplest common model---one can think of --- is linear, but  linear systems often fall short in accurately capturing the complex dynamics of nonlinear systems. In this work, we propose using quadratic systems as the common structure, inspired by the \emph{lifting principle}. 
According to this principle, smooth nonlinear systems can be expressed as quadratic systems in suitable coordinates without approximation errors. However, finding these coordinates solely from data is challenging. Here, we leverage deep learning to identify such lifted coordinates using only data, enabling a quadratic dynamical system to describe the system's dynamics. Additionally, we discuss the asymptotic stability of these quadratic dynamical systems.
We illustrate the approach using data collected from various numerical examples, demonstrating its superior performance with the existing well-known techniques. 
}

\novelty{
	\begin{itemize}
		\item We propose a data-driven approach to learn nonlinear dynamical models.
		\item We utilize the lifting-principle, allowing to write smooth nonlinear systems as quadratic systems in appropriate coordinates, which we refer to as quadratic embeddings.
		\item We introduce the discovery of quadratic embeddings using deep neural networks.
\end{itemize}
} 

\code{It will be made available on \url{https://gitlab.mpi-magdeburg.mpg.de/goyalp/quadembed_nonlineardyns}}
\maketitle

\section{Introduction}\label{sec:introduction}

Mathematical modeling is a crucial part of engineering design for control, optimization, and forecasting tasks. For those studies, it is important to have models that accurately represent dynamic behavior. Traditionally, such models are built based on fundamental principles or hypotheses or by field experts. While successful in many cases, it is challenging for complex processes like climate dynamics and modern robotics. Therefore, there is a growing need for data-driven modeling of dynamic systems using time-series data. Advances in sensor technology provide access to extensive data, making the development of such approaches more attractive and feasible.

Data-driven modeling of dynamical systems from time-series data has been an active research field; see, e.g., \cite{crutchfield1987equations,bongard2007automated,yao2007modeling,schmidt2011automated,rowley2009spectral,schmid2010dynamic,morBenGW15,morPehW16,rudy2017data,takeishi2017learning,morGoyB22a,morMayA07,morDrmGB15a,nakatsukasa2018aaa}. With the powerful approximation capability of deep neural networks (DNNs), several methods to model dynamic behavior using DNNs have been proposed \cite{rico1994continuous,rico1995nonlinear,gonzalez1998identification,mardt2018vampnets,vlachas2018data,champion2019data,lusch2018deep,chen2018neural,morGoyB21b}. One notable advantage of DNN-based modeling is the ability to approximate arbitrary functions without the need for manually designed features. 
Instead, DNNs learn essential features by identifying patterns in the data.
However, while neural networks excel at capturing complex dynamical behavior, applying them to engineering tasks like control and optimization can be challenging due to the \emph{black box} nature of DNN models, particularly when DNNs directly model the vector fields of dynamical systems, e.g., by using the NeuralODEs framework \cite{chen2018neural}.
Therefore, there is a need for simpler analytic dynamical models that can still capture complex dynamic behavior.

Linear models offer simplicity for various engineering studies, such as stability and control analysis and prediction, and extensive tools are also developed to carry out those studies for linear systems \cite{ogata2010modern,aastrom2021feedback,brunton2022data}. 
However, linear models often fail to capture intricate nonlinear dynamical behavior accurately. To address this limitation, additional rationales are required to incorporate nonlinear behavior. In this direction, the Koopman operator theory enables the representation of nonlinear systems as linear ones in an infinite-dimensional Hilbert space \cite{koopman1931hamiltonian}. Techniques like Dynamic mode decomposition (DMD) and its variants aim at approximating the infinite-dimensional Koopman operator in a finite-dimensional space using suitable coordinates or observables \cite{schmid2010dynamic,williams2015data,takeishi2017learning,morKutBBetal16,li2017extended,morBenHM18}. 
Despite the success of DMD in various applications, identifying appropriate observables for complex processes to approximate the Koopman operator remains challenging. 
An alternative approach involves using DNNs to discover a finite-dimensional approximation of the Koopman operator or linear embeddings \cite{lusch2018deep}. While promising, this method faces difficulties when dealing with models possessing continuous parts in their eigenspectrum. In such cases, determining a coordinate system and a finite-dimensional linear operator becomes difficult, see, e.g., \cite{lusch2018deep,mezic2020spectrum}. 

On the other hand, there exists a lifting principle that enables to write nonlinear dynamic models as quadratic models in lifted coordinates in \emph{finite-dimensions} \cite{savageau1987recasting,papachristodoulou2005analysis, morGu11}. This contrasts with Koopman operator theory, which might require an \emph{infinite-dimensional}  coordinate system. Determining such coordinates is feasible when we have analytical physics-based models for nonlinear dynamics, see, e.g., \cite{morGu11,morBenB15,qian2020lift}. However, our aim is to discover such coordinates in a data-driven framework since we only assume access to \emph{data}.

\begin{figure}[!tb]
	\includegraphics[width = 1\textwidth]{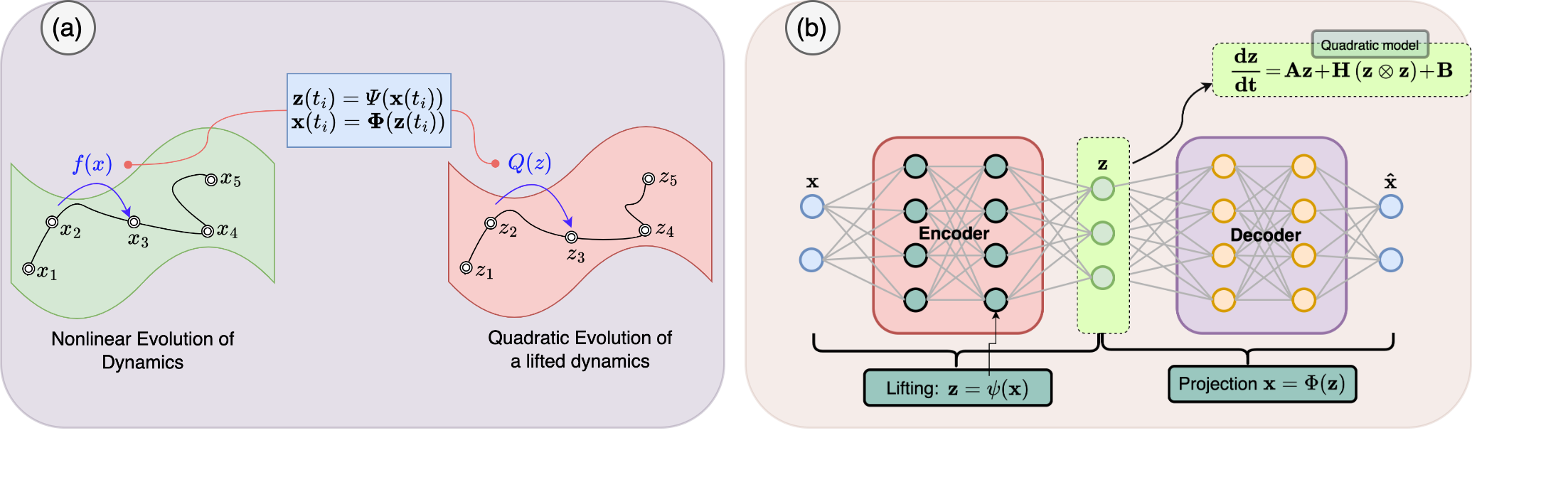}
	\caption{The figure illustrates that nonlinear dynamical systems can be written as quadratic dynamical systems in an appropriate finite-dimensional lifted coordinate system. In the right plot, we depict a neural network architecture to learn a  lifted coordinate system that has the desired quadratic embeddings.}
	\label{fig:method_pictorial}
\end{figure}

To that end, in this work, we present a DNN framework to learn the lifted coordinates using data. Our approach involves an autoencoder to discover lifted coordinates in which the dynamics  of the process generating the times series data can be governed by a quadratic dynamical model. In this paper, we refer to those lifted coordinates as \emph{quadratic embeddings}. We depict our objective in \Cref{fig:method_pictorial}. To achieve this, we simultaneously determine the parameters of the autoencoder and quadratic models for the lifted coordinates. Additionally, we delve into the asymptotic stability of quadratic models, explaining how the dynamics of embeddings are characterized and how stability can be ensured through appropriate parameterization. We demonstrate the success of our approach to learning dynamics via quadratic embeddings using two low-dimensional examples: the nonlinear pendulum and the Lotka-Volterra examples.

Additionally, we explore scenarios where data originates from complex, high-dimensional dynamic systems. The considered systems are known to evolve within a lower-dimensional subspace. When we have knowledge of physics-based models, it is possible to identify this subspace and associated reduced models using linear or nonlinear manifold projections; see, e.g., \cite{morSchVR08,lee2020model,barnett2022quadratic}. Another approach involves integrating the lifting principle  \cite{qian2020lift}, followed by constructing quadratic models in a reduced-dimensional space. Notably, the technique \cite{qian2020lift} employs a linear projection to determine a low-dimensional space, which is extended to quadratic manifold projection in \cite{geelen2023operator}. It is important to highlight that all these methods rely on physics-based discretized models or information at the partial differential equation (PDE) level. In our work, we discuss the construction of low-dimensional coordinates for high-dimensional systems using a convolution-based autoencoder so that a quadratic system can describe the dynamics of the low-dimensional coordinates. We illustrate this approach using the nonlinear Burgers example.

The remainder of the paper is structured as follows. The next section briefly recaps the lifting principle to determine quadratic embeddings for nonlinear systems. Section 3 presents a data-driven approach to learning lifted coordinates and the corresponding quadratic models. 
In Section 4, we discuss a parameterization for quadratic systems to ensure global asymptotic stability. Section 5 demonstrates the proposed method using various numerical examples. We finally conclude the paper with a summary and future avenues for further research in Section 6.

\section{Quadratic-embeddings for nonlinear systems}
McCormick, in 1976 \cite{mccormick1976computability}, proposed a convex relaxation to solve nonlinear non-convex optimization problems. The central idea lies in lifting a given nonlinear non-convex optimization problem into a high-dimensional space using auxiliary variables. Despite increasing the dimension of the original problem, we can solve the lifted problem more efficiently than the original one.
A similar philosophy has been developed in the context of dynamical systems. Casting nonlinear systems as polynomial systems can ease numerical analysis and control design tasks, see, e.g., \cite{savageau1987recasting,papachristodoulou2005analysis}. Furthermore, for constructing low-dimensional models for large-scale nonlinear dynamical systems, Gu \cite{morGu11} proposed a methodology that first rewrites nonlinear systems as quadratic systems and then employs reduced-order modeling techniques for quadratic systems. In addition, once we have a quadratic representation of a nonlinear system, we can apply design tools, such as control design, developed for quadratic systems (see, e.g., \cite{amato2007state,tognetti2021output}) to nonlinear systems as well. 
For given nonlinear analytic systems,  there are approaches to rewrite them as polynomial or quadratic systems \cite{savageau1987recasting,papachristodoulou2005analysis, morGu11}. There it, has been shown that nonlinear systems consisting of basic elementary smooth nonlinear functions (e.g., trigonometric, exponential, rational) or a composition of these elementary functions can always be rewritten as quadratic systems in lifted coordinates without any approximation error. Moreover, the dimension of the lifted coordinates increases linearly with respect to the dimension of the original nonlinear system, and one can also determine bounds on the dimension of the lifted coordinate system. 

Inspired by the lifting philosophy, we can transform all smooth nonlinear systems satisfying the above assumptions as a quadratic system. Furthermore, we can construct the original state vector or measurement space using a linear projection of the lifted variables. To illustrate the approach, we consider a simple rational nonlinear system often appearing in biological applications. The nonlinear model is as follows:
\begin{equation}
	\dot{\bx}(t) = -\dfrac{\bx(t)}{1+\bx(t)}.
\end{equation}
To write the above nonlinear system as a quadratic system, we define lifted coordinates as follows:
\begin{equation}
	\cL(\bx) := \begin{bmatrix}
		\bx \\[5pt] \dfrac{1}{1+\bx} \\[10pt] \dfrac{\bx}{(1+\bx)^2}
	\end{bmatrix}\equiv \by.
\end{equation}
Note that $\bx$ can be constructed using a linear projection of $\by$; precisely, using a matrix $\bC = \left[1, 0, 0\right]$, i.e., $\bx = \bC\by$. Moreover, the differential equation for $\by$ can be given as 
\begin{equation}
	\begin{bmatrix} y_1 \\y_2\\y_3  \end{bmatrix} = \begin{bmatrix}
		-y_1y_2 \\ y_2y_3 \\ y_3(y_2-2y_3)
	\end{bmatrix},
\end{equation}
which is a quadratic system with  $y_i$ being the $i$-th component of $\by$. 
These discussions indicate that given an analytical form of a nonlinear system, we can find a lifted coordinate system to write the dynamics as quadratic systems. We highlight that the lifted variables are not unique. 
Since here we only assume to have access to data, we do not have the possibility of hand-designing  lifted coordinates. Even if we had hand-designed coordinates, they are not straightforward for complex processes. Therefore, we next discuss a data-driven approach to learning a suitable lifted coordinate system.
\section{Data-driven discovery of quadratic embeddings using deep learning}
The lifting principle suggests that smooth nonlinear systems can be rewritten as quadratic systems in a suitable lifted coordinate system. Additionally, this coordinate system allows the construction of the quantity of interest by a nonlinear projection, as depicted in \Cref{fig:method_pictorial}(a). Inspired by these observations, we formalize our problem in the following.

\begin{problem}\label{prob:lifting}
	Given data $\{\bx(t_1),\ldots, \bx(t_\cN)\}$ and derivative information $\{\dot\bx(t_1),\ldots, \dot\bx(t_\cN)\}$, we seek to determine a lifting $\bx \mapsto \bz $ such that 
	\begin{enumerate}[label=(\alph*)]
		\item $\bz$ satisfies
		\begin{equation}\label{eq:lifted_quadmodel}
			\dot{\bz}(t) = \bA\bz(t) + \bH\left(\bz(t)\otimes \bz(t)\right) + \bB,\quad\text{and}
		\end{equation}
		\item $\bx(t)$ can be recovered by a (potentially) nonlinear projection of $\bz(t)$, i.e., $\bx(t) = \Phi(\bz(t))$,
	\end{enumerate}
	for $t \in \{t_1,\ldots, t_\cN\}$.
\end{problem}
Since we only have access to data and have no information about the underlying physics-based model, we cannot identify the required lifted coordinates analytically. Therefore, we aim to learn suitable coordinates using \emph{only} data. 
For this, we utilize the profound approximation capabilities of DNNs to learn the lifting. To that end, we propose a particular autoencoder design, where the encoder is defined by a DNN denoted by $\Psi$, which is parameterized by $\theta$, such that it provides lifted coordinates for a given $\bx$, i.e., $\bz =\Psi(\bx)$. Moreover, the decoder is given by a nonlinear function that maps the lifted variables back to the space of measurements by a nonlinear function, i.e., $\bx = \Phi(\bz)$. The nonlinear function $\Phi(\cdot)$ is also constructed using a DNN, and its parameters are denoted by $\phi$.
Suppose we were to train the autoencoder alone. In that case, we may obtain any arbitrary mapping, and it is not necessary that the dynamics of the learned lifted coordinates can be described by a quadratic model. Therefore, it is essential to learn the parameters of the autoencoder to determine lifted coordinates so that among infinitely many possibilities, we can identify ones that fulfill both conditions (a) and (b) in \Cref{prob:lifting}.

Towards achieving our goals, we assume to have access to the derivative information of $\bx$. Then, we can compute the derivative of the lifted coordinates $\bz$ using the chain rule and automatic differentiation, i.e., $\dot{\bz}(t) = \nabla_\bx\Psi(\bx(t))\dot\bx(t)$, where $\nabla_\bx$ denotes the Jacobian with respect to $\bx$. We further enforce that a quadratic equation as in \eqref{eq:lifted_quadmodel} can also provide the derivative information of $\bz$; hence, we add the following term in the loss function:
\begin{equation}\label{eq:loss1}
	\cL_{\dot\bz\dot\bx} = \dfrac{1}{\cN} \sum_{i = 1}^{\cN}\left\| \nabla_\bx\Psi(\bx(t_i))\dot\bx(t_i) -\left( \bA\bz(t_i) + \bH\left(\bz(t_i)\otimes \bz(t_i)\right) + \bB\right)\right\| \quad \text{with}~~\bz(t_i) = \Psi(\bx(t_i)).
\end{equation}
Moreover, as in \cite{champion2019data}, we can also reconstruct the derivative information of the original variable $\bx$ using the derivative information of the lifted variable $\bz$. Using the nonlinear mapping between $\bx$ and $\bz$ via the decoder, we have 
\begin{equation}\label{eq:loss2}
	\dot\bx = \nabla_\bz \Phi(\bz) \dot\bz =  \nabla_\bz \Phi(\bz)\left(\bA\bz + \bH\left(\bz\otimes\bz\right)+\bB\right).
\end{equation}
This gives rise to the second part of the loss function as follows:
\begin{equation}
	\cL_{\dot\bx\dot\bz} = \dfrac{1}{\cN} \sum_{i=1}^{\cN} \|	\dot\bx(t_i) -  \nabla_\bz \Phi(\bz(t_i))\left(\bA\bz(t_i) + \bH\left(\bz(t_i)\otimes\bz(t_i)\right)+\bB\right)\|.
\end{equation}
Naturally, the reconstruction of the original state $\bx$ from $\bz$ can be obtained as follows:
\begin{equation}
	\bx = \Phi(\Psi(\bx)),
\end{equation}
yielding the third element of the loss function 
\begin{equation}\label{eq:loss3}
	\cL_{\text{encdec}}	= \dfrac{1}{\cN} \sum_{i=1}^{\cN}  \|\bx(t_i) - \Phi(\Psi(\bx(t_i)))\|.
\end{equation}
Combining all these elements \cref{eq:loss1,eq:loss2,eq:loss3}, we have the following total loss:
\begin{equation}\label{eq:total_loss}
	\cL = \lambda_1\cL_{\text{encdec}} +\lambda_2 \cL_{\dot\bx\dot\bz}  +\lambda_3 \cL_{\dot\bz\dot\bx} ,
\end{equation}
where $\lambda_{\{1,2,3\}}$ are the hyper-parameters. Moreover, the norms in \cref{eq:loss1,eq:loss2,eq:loss3} are a weighted sum of the Frobenius norm and $l_1$-norm, i.e., $\|\cdot\| = 0.5 \|\cdot \|_F + 0.5 \|\cdot \|_{l_1}$.
We optimize simultaneously the parameters of the autoencoder  $\{\ensuremath{\mathbf{\theta}},\mathbf{\phi}\}$ and the system matrices $\{\bA,\bH,\bB\}$, minimizing the loss $\cL$ in \eqref{eq:total_loss}.
Furthermore, we highlight that once the autoencoder is trained and system matrices defining the dynamics are obtained, we query the nonlinear encoder only to get the corresponding initial condition for the lifted variable $\bz$. We can obtain the trajectory of $\bz$ using a desired numerical integration method, so the trajectory of $\bx$ can be determined using the decoder. Hence, the acquired dynamic modeling is parsimonious by construction compared to DNNs, which directly model the vector field of $\bx$. 

\section{Asymptotic stability-guaranteeing quadratic embeddings}\label{sec:stability}
In the previous section, we discussed the problem of quadratic embeddings for nonlinear systems. This section discusses an approach to ensuring the asymptotic stability property of the quadratic embeddings. Precisely, we are interested in guaranteeing the asymptotic stability of the quadratic system that describes the dynamics of the identified embeddings. Such property might be essential to have, particularly when the original systems from which data come are asymptotically stable. 

To achieve our stability goals, we use the results presented in \cite{goyal2023guaranteed}. The work proposes a parameterization of quadratic systems, which ensure global asymptotic stability, \emph{not just local}, of quadratic systems. The results are summarized in the following lemma. 
\begin{lemma}[\cite{goyal2023guaranteed}]\label{lemma:stability}
	Consider a quadratic system as follows:
	\begin{equation}
		\dot{\by}(t) = \bA\by + \bH(\by\otimes \by),
	\end{equation}
	where $\bA \in \Rnn$ and $\bH \in \R^{n\times n^2}$. Assume that $\bA$ can be written as $\bJ -\bR$ with $\bJ = -\bJ^\top$ and $\bR = \bR^\top \succ 0$, and $\bH = \begin{bmatrix}	\bH_1, \ldots, \bH_n \end{bmatrix}$ with $\bH_i = -\bH_i ^\top$, $i \in \{1,\ldots, n\}$. Then, the quadratic system is asymptotically stable, i.e., $\lim_{t\rightarrow \infty}\by(t) \rightarrow 0$. Furthermore, if $\bR = \bR^\top \succeq 0$, then $\|x(t)\|_2 \leq \|x_0\|_2$, where $\by(t)$ is the solution at time $t$ for a given initial condition $\by_0$.
\end{lemma}
We utilize the results from \Cref{lemma:stability} and combine them with the approach presented in the previous section to identify quadratic embeddings. Hence, we modify \Cref{prob:lifting} as follows:
\begin{problem}\label{prob:lifting_stability}
	Given data $\{\bx(t_1),\ldots, \bx(t_\cN)\}$ and derivative information $\{\dot\bx(t_1),\ldots, \dot\bx(t_\cN)\}$, we seek to determine a lifting $\bx \mapsto \bz $ such that 
	\begin{enumerate}[label=(\alph*)]
		\item $\bz$ satisfies
		\begin{equation}\label{eq:lifted_quadmodel_stable}
			\dot{\bz}(t) = \bA\bz + \bH\left(\bz\otimes \bz\right),
		\end{equation}
		where $\bA = \bJ -\bR$ with $\bJ = -\bJ^\top$ and $\bR = \bR^\top \succ 0$, and $\bH = \begin{bmatrix}	\bH_1, \ldots, \bH_n \end{bmatrix}$ with $\bH_i = -\bH_i ^\top$, $i \in \{1,\ldots, n\}$. 
		\item $\bx$ can be recovered using $\bz$, i.e., $\bx = \Phi(\bz)$. 
	\end{enumerate}
\end{problem}
When quadratic embeddings are determined so that it solves \Cref{prob:lifting_stability}, we can guarantee the stability of the embedded dynamics by construction. The rest of the steps for the construction of losses and the usage of the autoencoder remain the same as in the previous section. 

\begin{remark}
	If we set $\bR = 0$ in \Cref{prob:lifting_stability}, then it can be easily shown that the energy, defined by the $2$-norm of the vector, is conserved. This means that $\dfrac{d}{dt} \|z(t)\|_2 = 0$. This is interesting when learning dynamical systems that are energy-preserving, e.g., the Schr\"odinger equation or friction-less multi-body dynamics. 
\end{remark}

\section{Numerical Results}

We demonstrate our approach with three examples: a nonlinear pendulum, a dissipative Lotka-Volterra system, and a 1D Burgers' model. All the information about data generation, training, autoencoder design, and hyper-parameters is given in \ref{sec:appendix}. We refer to our proposed methodology as \quadembs. As for all three models, we expect globally stable behavior, we enforce the asymptotic stability of the embeddings via parameterization as discussed in \Cref{sec:stability}.
We compare our approach with two existing methodologies. The first one is proposed in \cite{lusch2018deep}, inspired by the Koopman theory, which aims to learn universal linear embeddings for nonlinear dynamics. Note that \cite{lusch2018deep} discusses a universal linear embedding approach for discrete systems, which we slightly adopt for the continuous setting. We denote this approach as \linearembs. Additionally, \linearembs~also requires an autoencoder structure; thus, we design it as for \quadembs. The second approach, we consider, is the operator inference approach \cite{morPehW16}, which aims to learn polynomial systems given data. We restrict the degree of the polynomial to two. We refer to this approach as \quadopInf. 

\subsection{A nonlinear pendulum} In the first example, we consider a nonlinear damped pendulum, which is governed by the following second-order equation:
\begin{equation}
	\ddot{\bx}(t)= -\sin(\bx(t)) - 0.025\dot\bx,
\end{equation}
which, in the first-order companion form, can be written as 
\begin{equation}
	\begin{bmatrix} \dot{\bx}_1(t) \\ \dot{\bx}_2(t) \end{bmatrix}= \begin{bmatrix} -\sin(\bx_2(t)) - 0.025\bx_1(t) \\ \bx_1(t) \end{bmatrix},
\end{equation}
where $\bx_1(t) = \dot\bx(t)$ and $\bx_2(t) = \bx(t)$. We collect the data in the time interval $t \in \left[0,25\right]$s with $50$ different initial conditions. For each initial condition, we take $100$ equidistant sample points in the given time interval. We randomly choose an initial position and velocity of the pendulum in the range of $[-3,3]$. We assume to know the gradient of $\bx_1$ and $\bx_2$. To employ \quadembs, we design an autoencoder with the hyper-parameters given in \Cref{tab:hyperparameters_lowdimensional} and set the dimension of the lifted variables to three. The same autoencoder design is used for \linearembs~and linear embeddings of the latent dimension three are learned. Additionally, \quadopInf~does not involve an autoencoder, and it aims to learn the dynamics in the original coordinate systems. Hence, the dynamical model learned using \quadopInf~is of order two. 

\begin{figure}[tb]
	\centering
	\includegraphics[width = 0.495\linewidth]{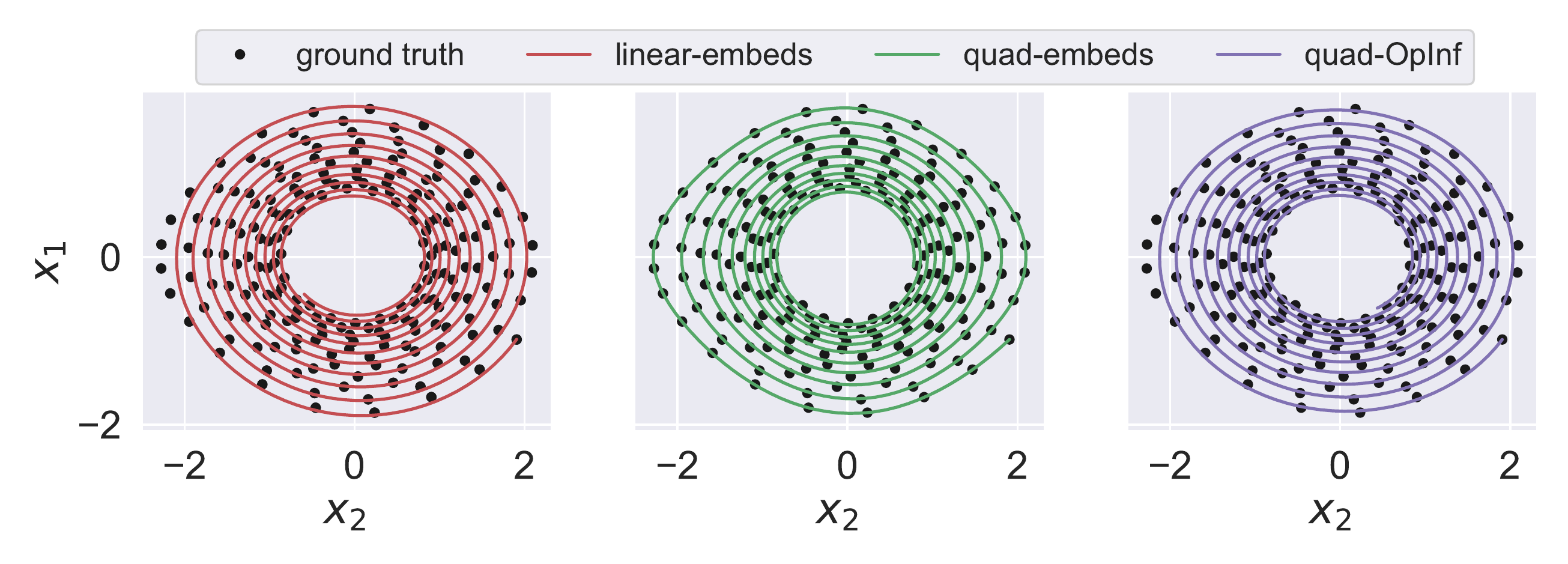}
	\includegraphics[width = 0.495\linewidth]{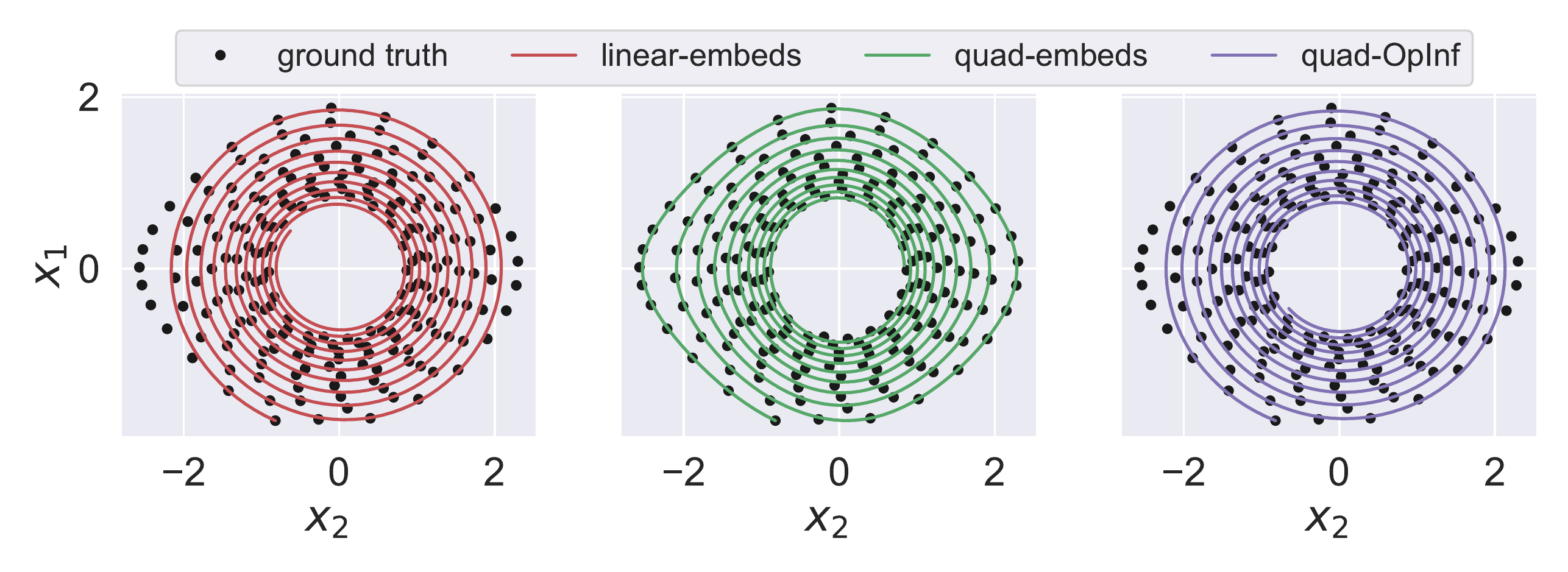}\\
	\includegraphics[width = 0.495\linewidth]{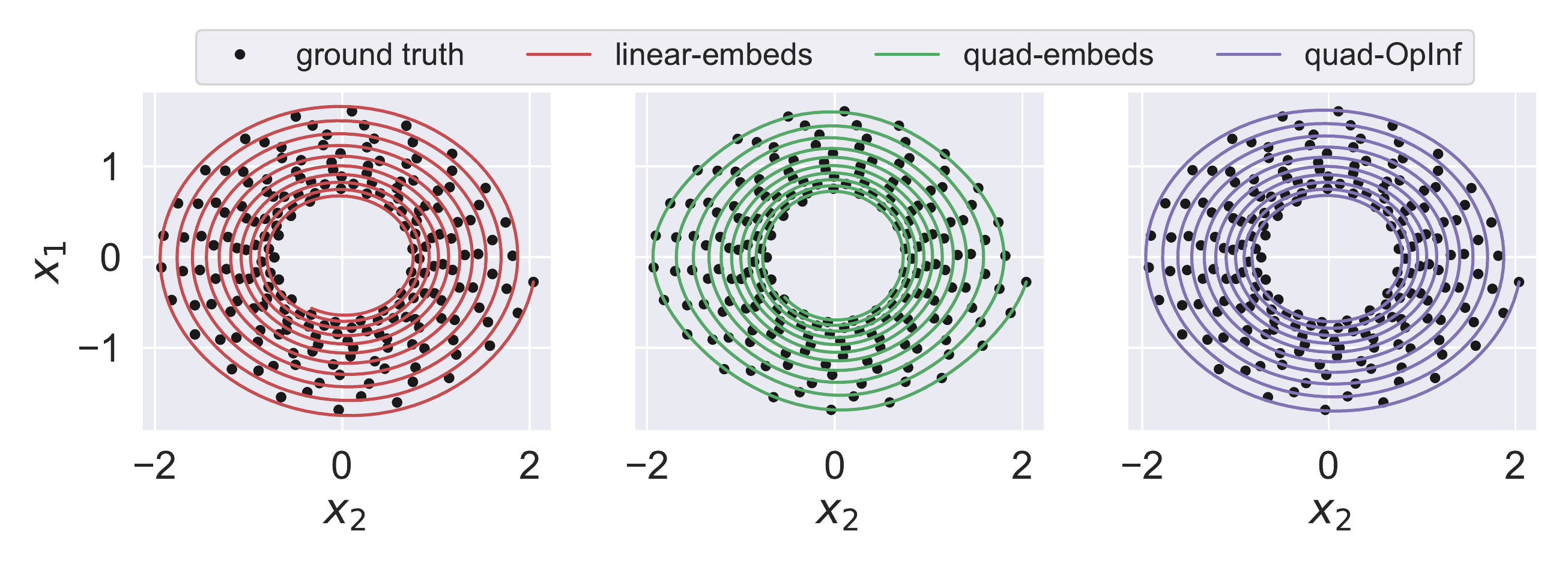}
	\includegraphics[width = 0.495\linewidth]{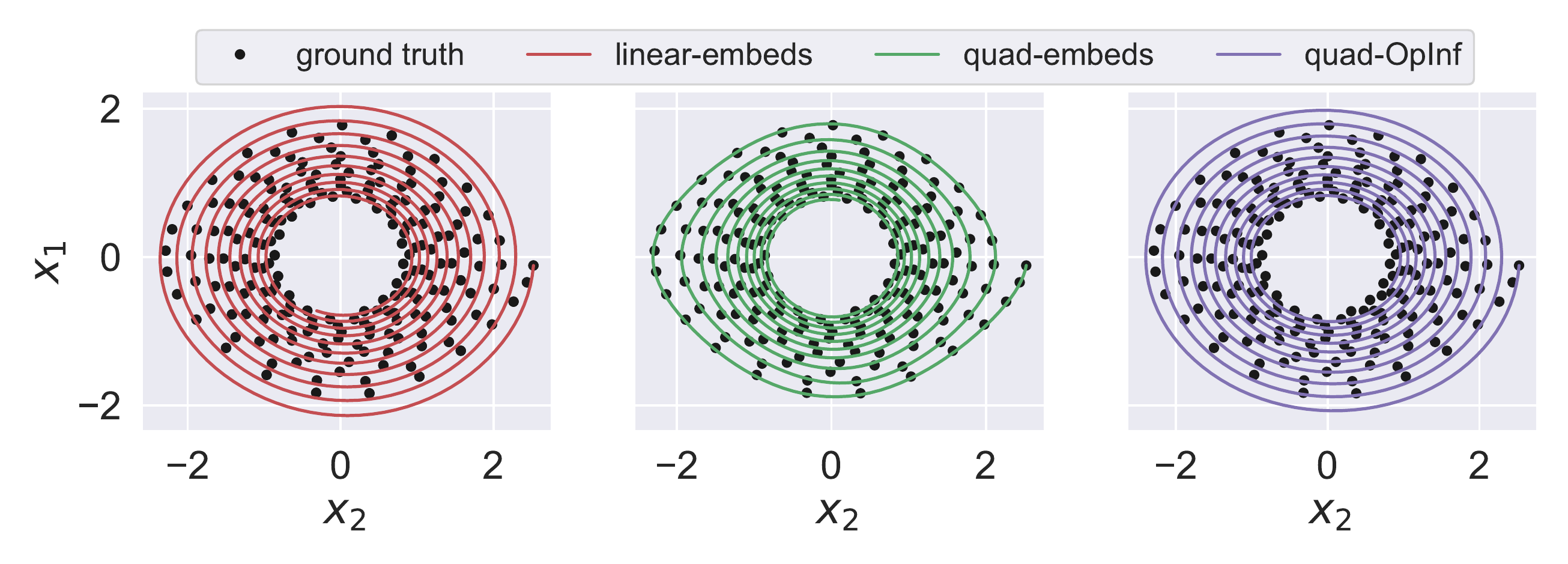}
	\caption{Nonlinear pendulum example:  A comparison of the trajectories obtained using \linearembs, \quadembs, and \quadopInf~methods with the ground truth ones on the testing data is presented.  }
	\label{fig:pendulum_example}
\end{figure}

Having learned models using these different methods, we test their performances using $100$ initial conditions, which have not been part of the training. Moreover, the time interval for testing is considered to be $t \in [0,75]s$, which is three times longer than the training one, and we take $2,000$ equidistant sample points in the testing interval. Next, in \Cref{fig:pendulum_example}, we show comparisons for four trajectories, which are chosen for which \quadembs~performs the worst. Note that performance is measured based on the following criteria:
\begin{equation}\label{eq:error_measure}
	\cE(\bx_0) = \log_{10} \left(\texttt{median}\left( \left(\bX^{(\bx_0)}_{\texttt{ground-truth}} - \bX^{(\bx_0)}_{\texttt{learned}}   \right) \right)^2 \right),
\end{equation}
where $\bX^{(\bx_0)}_{\texttt{ground-truth}}$ and  $\bX^{(\bx_0)}_{\texttt{learned}}$, respectively, contain solutions using the ground truth model and the learned models for a given test initial condition ${\bx_0}$. We notice that \linearembs~and \quadopInf~poorly perform to compare the dynamics. A potential reason for the failure of \linearembs~might be associated with the continuous spectrum of the pendulum system. Furthermore, as \quadopInf~learns a quadratic model in the original measurement coordinates without any transformation, the quadratic model is not capable of capturing the dynamics completely. On the other hand, the proposed methodology \quadembs~learned the dynamics of the original systems accurately by learning suitable quadratic embeddings using an autoencoder, thus illustrating the power of learning suitable quadratic embeddings for nonlinear systems with continuous spectra. 

\begin{figure}[tb]
	\includegraphics[width = 0.45\linewidth]{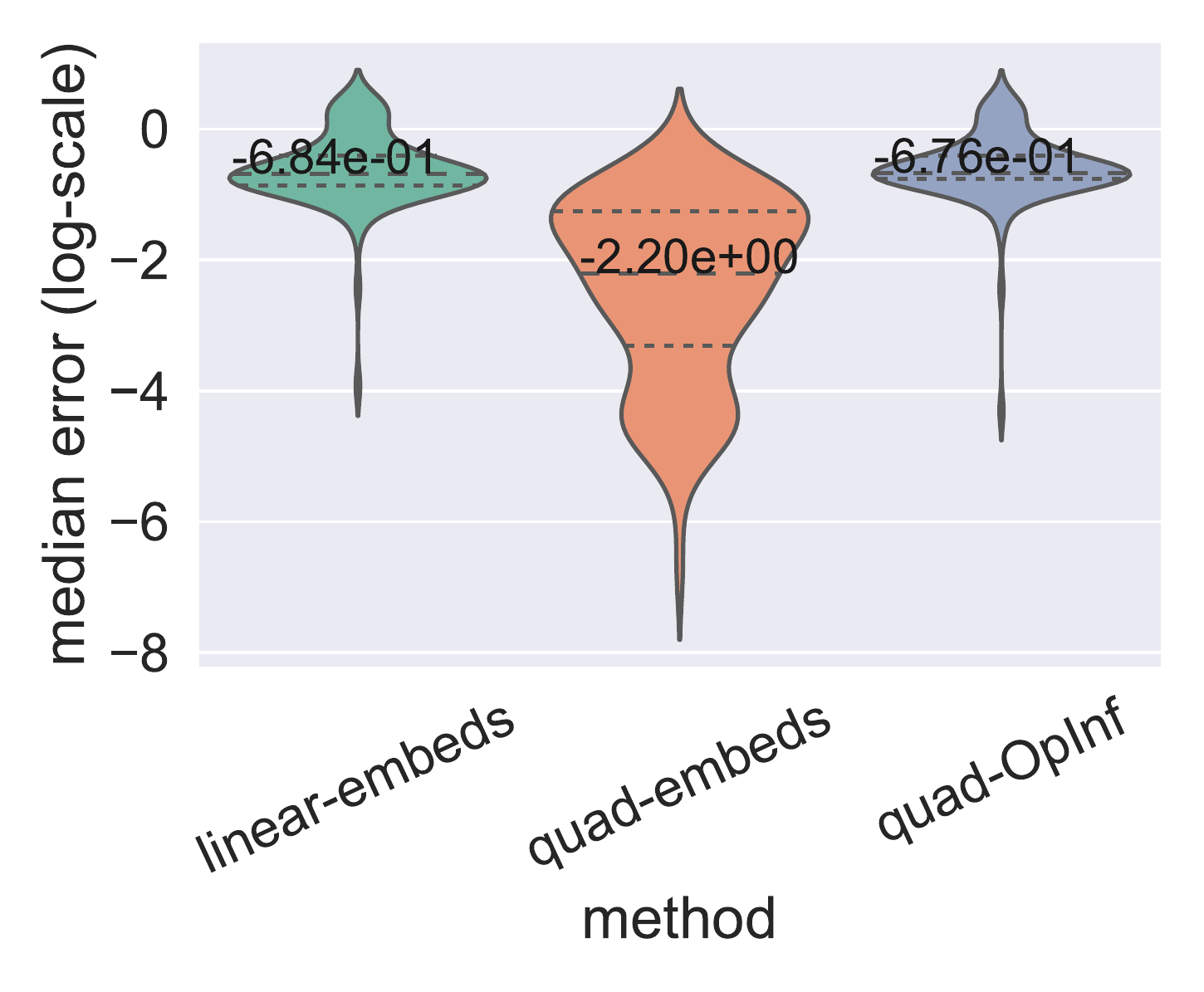}
	\caption{Nonlinear pendulum example:  The figure shows a qualitative comparison of the performance of  \linearembs, \quadembs, and \quadopInf~on the testing data based on the measure \eqref{eq:error_measure}. Note that the error measure \eqref{eq:error_measure} contains the \emph{log}; hence, more negative is the values, better the method perform better.}
	\label{fig:pend_error_measure}
\end{figure}

Furthermore, for a detailed performance comparison of all three methods, we compute the error using  \eqref{eq:error_measure} for each test initial condition and show them using the violin-plot in \Cref{fig:pend_error_measure}. The figure clearly indicates a superior performance of the proposed methodology. 

\subsection{Dissipative Lotka-Volterra example}
In our second example, we consider the Lotka-Volterra example. The example is often used as a benchmark for a Hamiltonian system \cite{nutku1990hamiltonian}; however, we, here, consider its dissipative form, which is governed using the  following equations:
\begin{equation}
	\begin{bmatrix} \dot \bq(t) \\ \dot \bp(t) \end{bmatrix} =  \begin{bmatrix}  - e^{\bp} - 0.05\cdot \bq + 1 \\  e^{\bq} - 0.05\cdot \bp -2  \end{bmatrix},
\end{equation}
where $\bq$ and $\bp$ represent position and momentum quantities, respectively. We collect the data in the time interval $t \in \left[0,10\right]$s with \emph{only} $10$ different initial conditions. We take $200$ equidistant sample points for each initial condition in the considered time interval. We randomly choose an initial position and momentum in the range of $[-1.5,1.5]$. Furthermore, analogous to the previous example,  we assume to have the gradient information for $\bq$ and $\bp$. To employ \quadembs, we design an autoencoder with the parameters given in \Cref{tab:hyperparameters_lowdimensional} and set the dimension of the lifted coordinate system to three. The same autoencoder design is used for \linearembs~and a linear operator is learned with the dimensional of the latent representation being three. Since \quadopInf~does not involve an autoencoder and it learns the dynamics using the measured state coordinates, the dynamical model learned using \quadopInf~is of order two. 

\begin{figure}[tb]
	\centering
	\includegraphics[width = 0.495\linewidth]{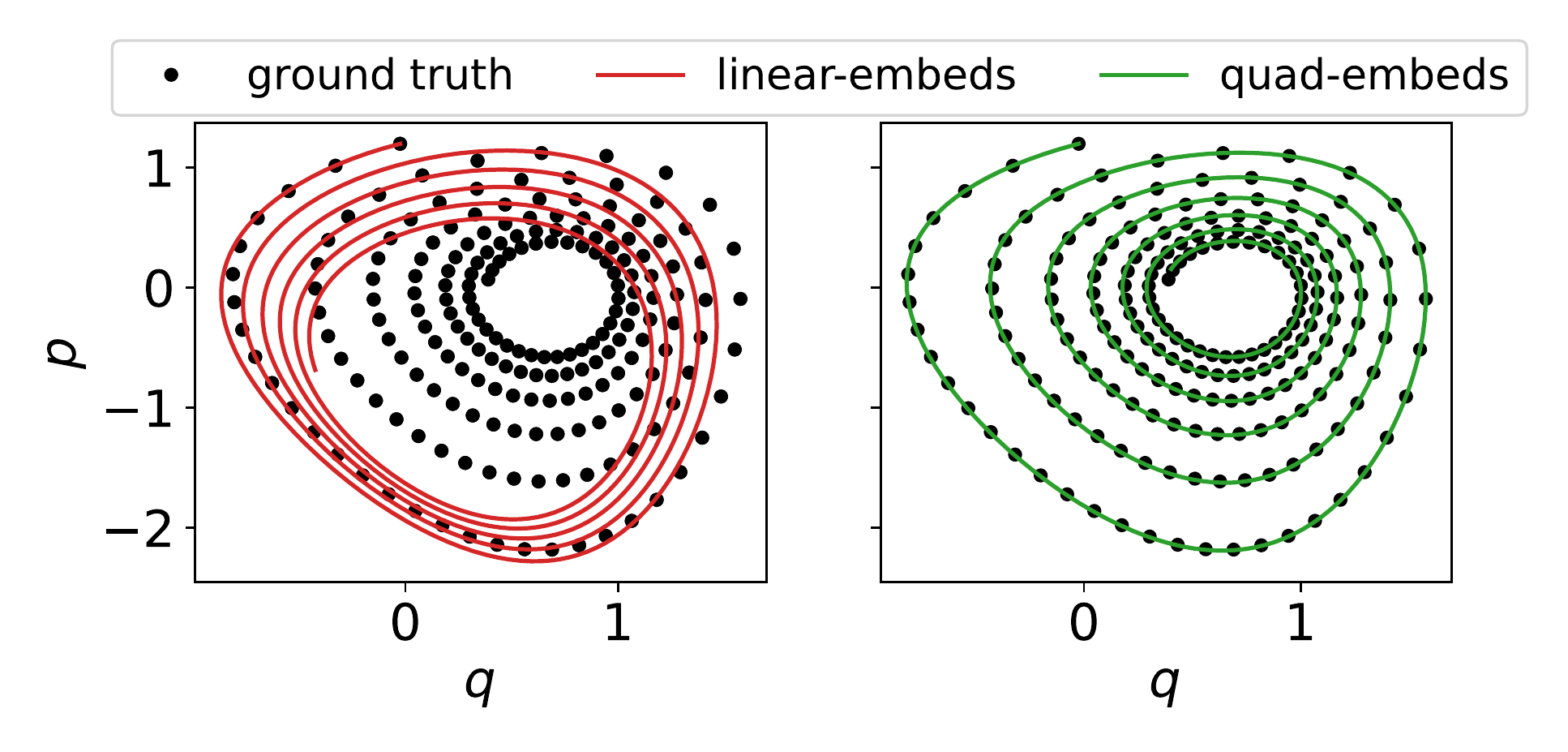}
	\includegraphics[width = 0.495\linewidth]{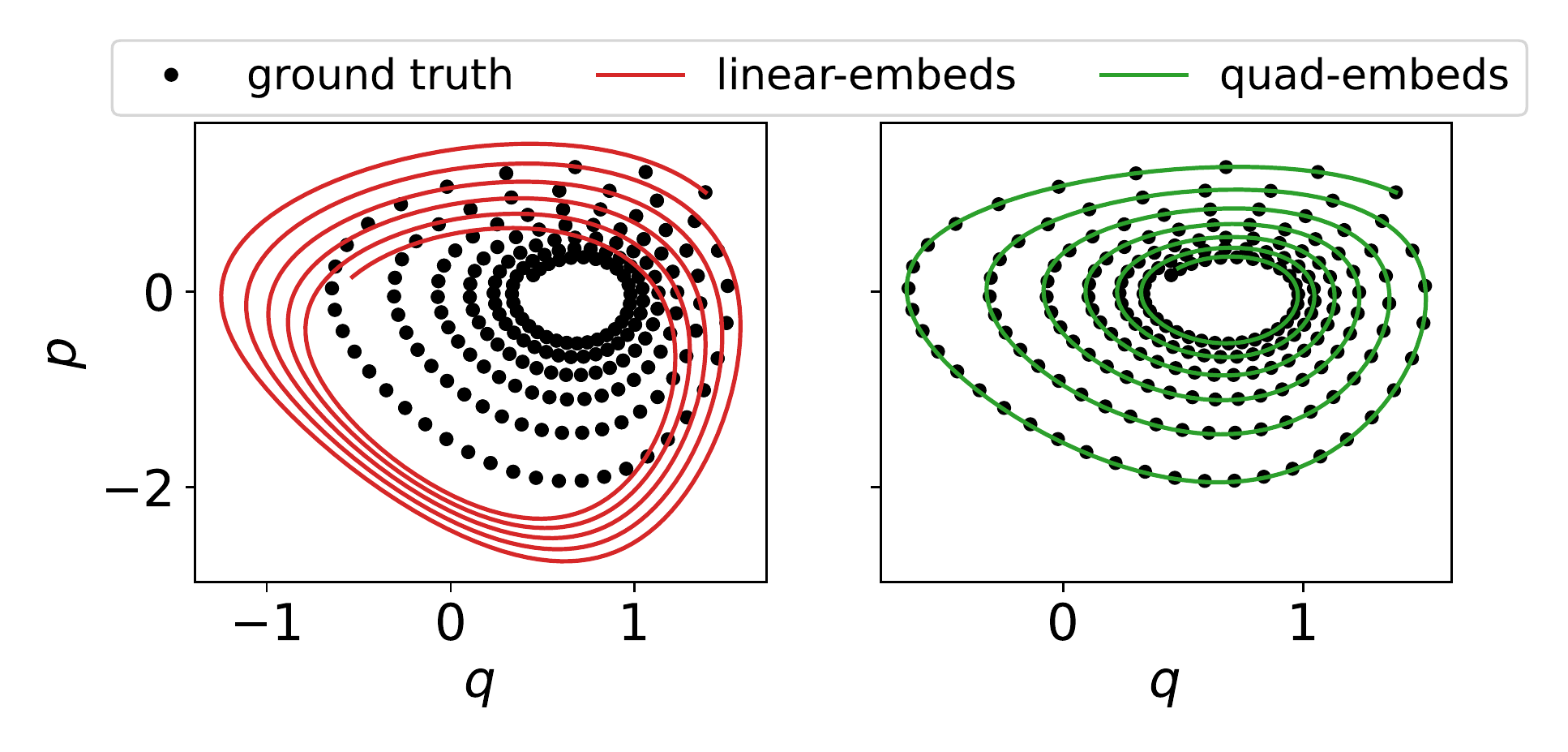}\\
	\includegraphics[width = 0.495\linewidth]{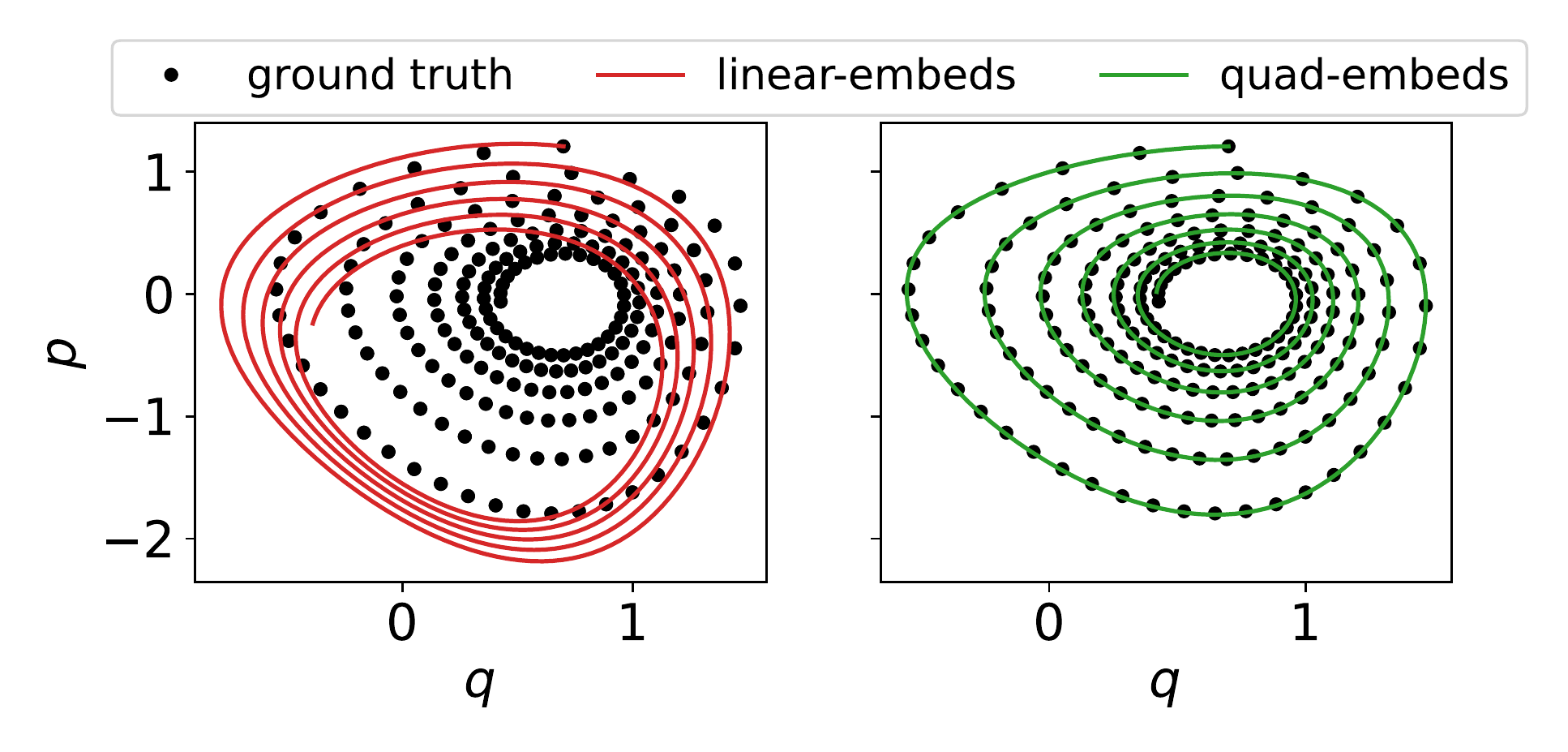}
	\includegraphics[width = 0.495\linewidth]{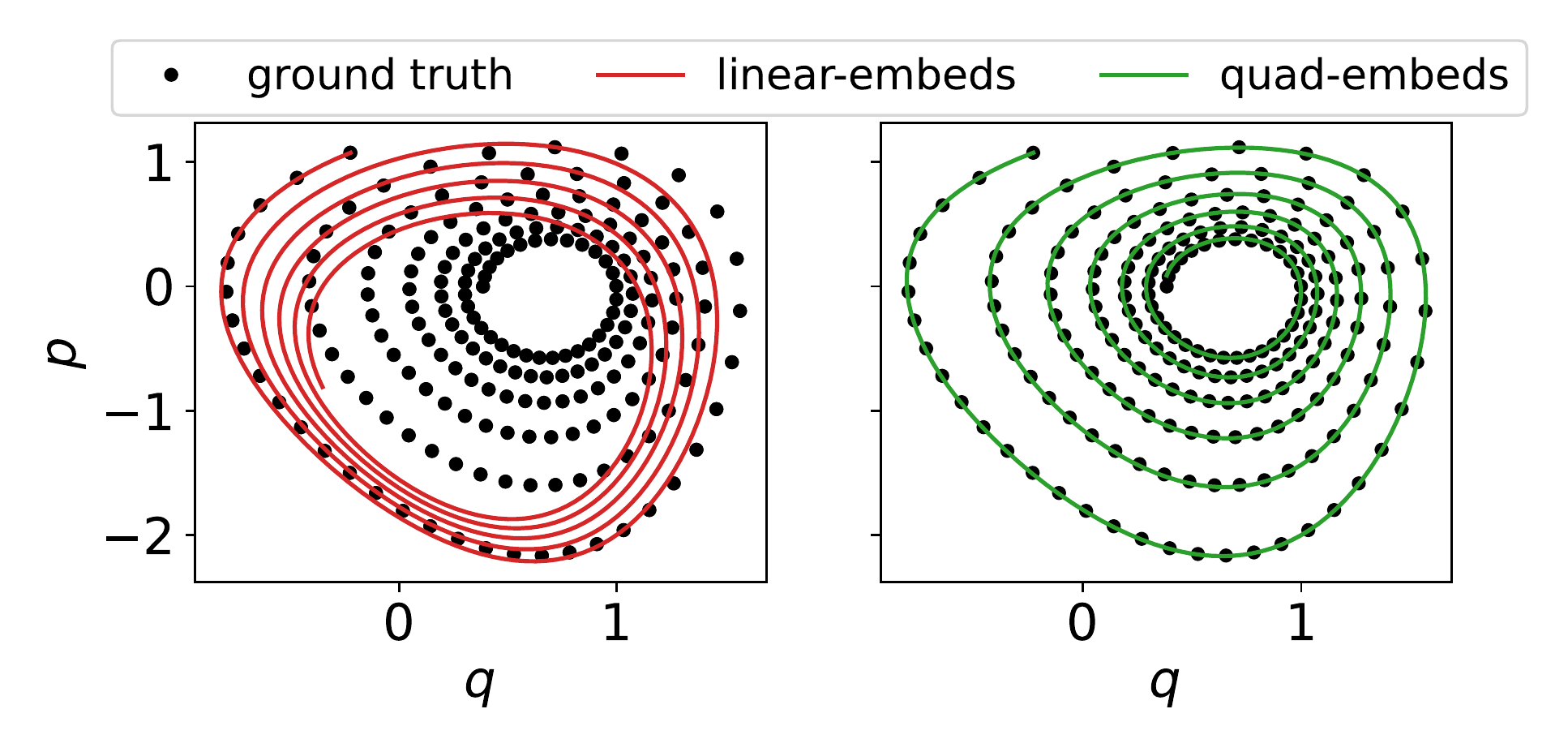}
	\caption{Lotka-Volterra example:  A comparison of the trajectories obtained using \linearembs, \quadembs, and \quadopInf~methods with the ground truth ones on the testing data is presented. Note that \quadopInf~yields unstable trajectories; therefore, they are not shown in the figure.}
	\label{fig:LV_example}
\end{figure}

Having learned models using these different methods, we test their performances using $100$ initial conditions, which have not been part of the training. Moreover, the time interval for testing is considered to be $t \in [0,30]$s, which is three times longer than the training one, and we take $4,000$ equidistant sample points in the testing time interval. Similar to the previous example, in \Cref{fig:LV_example}, we show comparisons for four trajectories, which are chosen where \quadembs~performs the worst based on the criterion in \eqref{eq:error_measure}. We notice that \linearembs~and \quadopInf\ completely fail to capture the dynamics. Potential reasons for the failures of \linearembs~and \quadopInf~could be the same as for the previous example. Furthermore, we highlight that \quadopInf~yields unstable trajectories; therefore, the plots are not shown for \quadopInf~in \Cref{fig:LV_example}. Additionally, we note that  $50$ trajectories out of $100$ testing ones are unstable for \quadopInf. On the other hand, the proposed methodology \quadembs~learned the dynamics of the original systems accurately by learning suitable quadratic embeddings using an autoencoder, and all the trajectories are stable by construction. 

\begin{figure}[tb]
	\includegraphics[width = 0.45\linewidth]{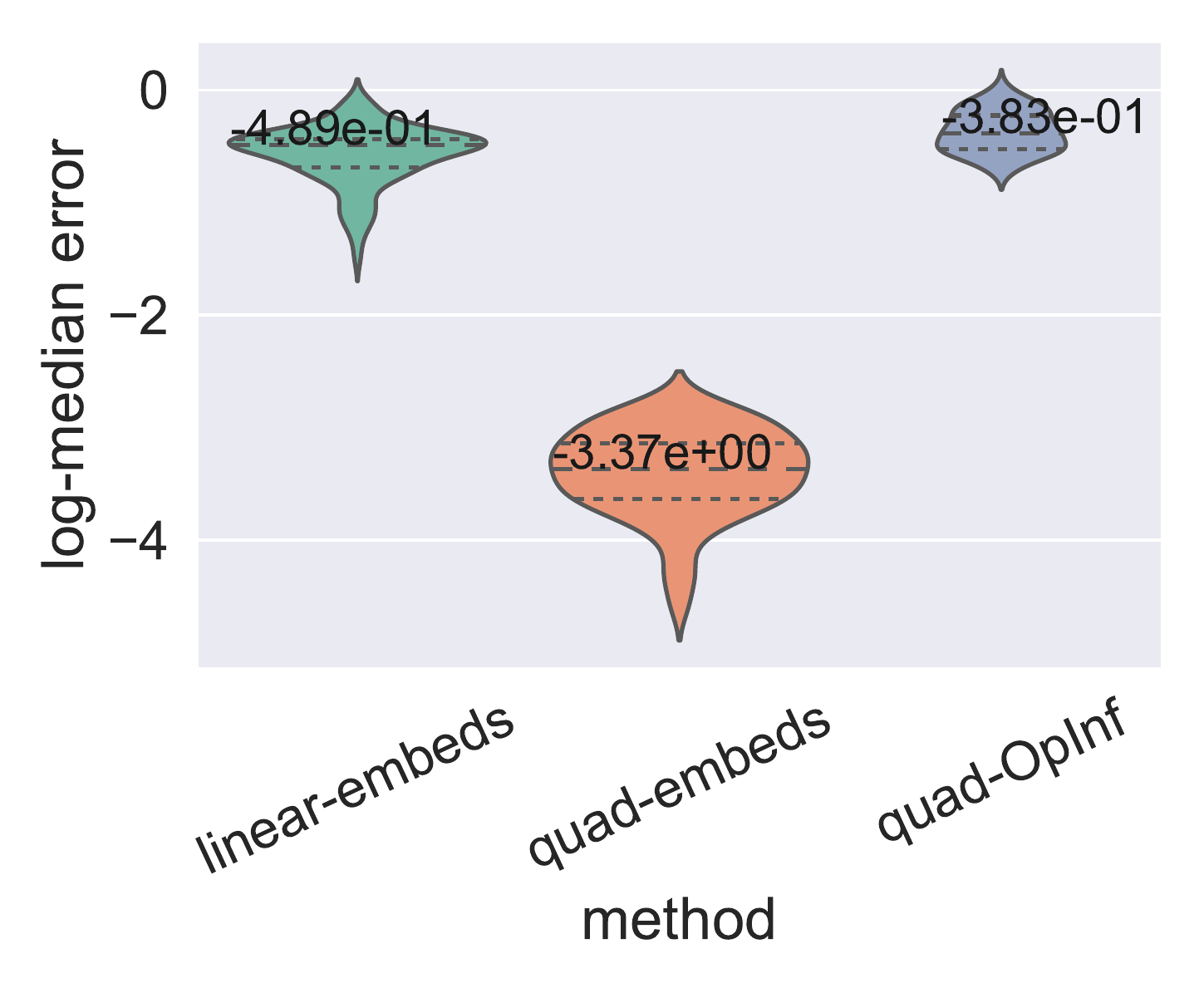}
	\caption{Lotka-Volterra example:  The figure shows a qualitative comparison of  \linearembs, \quadembs, and \quadopInf~on the testing data.  Note that while plotting the error for  \quadopInf, unstable trajectories are removed, which are as many as $50$ out of $100$ testing initial conditions.}
	\label{fig:LV_violinplot}
\end{figure}

Furthermore, for a detailed qualitative comparison of all three methods, we compute the error using \eqref{eq:error_measure} and show them using the violin-plot in \Cref{fig:LV_violinplot}. The figure clearly indicates a superior performance of the proposed methodology. 

\subsection{A High-dimensional Example: Nonlinear Burgers' Equations}
In our last example, we consider a somewhat higher example inspired by Burgers' equations. It is a one-dimensional PDE with the governing equations as follows:
\begin{subequations}
	\begin{align}
		&u_t + uu_{x} + u^3u_x = u_{xx}, \quad \text{with}~~ x \in (0,1) ~\text{and}~ t \in (0,T),\label{eq:burgers_governingequation}\\
		& u(0,\cdot) = 0, \quad \text{and}~ u(1,\cdot) = 0,\\
		& u(x, 0) = 10\cdot \sin(\pi x \cdot f) x (1-x), \label{eq:burgers_inits_conditions}
	\end{align}
\end{subequations}
where  $u_t, u_x$ and $u_{xx}$ represent the derivative of $u$ with respect to time $t$, the derivative of $u$ with respective to space $x$, and the double derivative of $u$ with respect to $u$, respectively;   $f \in \R$ is a frequency. We highlight the additional term $u^3u_x$ \eqref{eq:burgers_governingequation}, which makes the modified equations have quartic polynomial terms. We collect the data by considering various initial conditions, and for this, we vary the parameter $f$ in \eqref{eq:burgers_inits_conditions}. We take $13$ different values of $f$ equidistantly in the range $[2,3]$. Assuming the values of $f$ are sorted in increasing order, we consider th $3$rd, $6$th, $9$th, and $12$th values for the testing and the remaining nine values for training. Moreover, we discretize the PDE using a finite-difference scheme by considering $256$ points in space, and we integrate the discretized system in the time interval $[0,1.5]$s and collect $1001$ data points in the interval. 

\begin{figure}[tb]
	\centering
	\includegraphics[width = 0.75\linewidth,angle = 0]{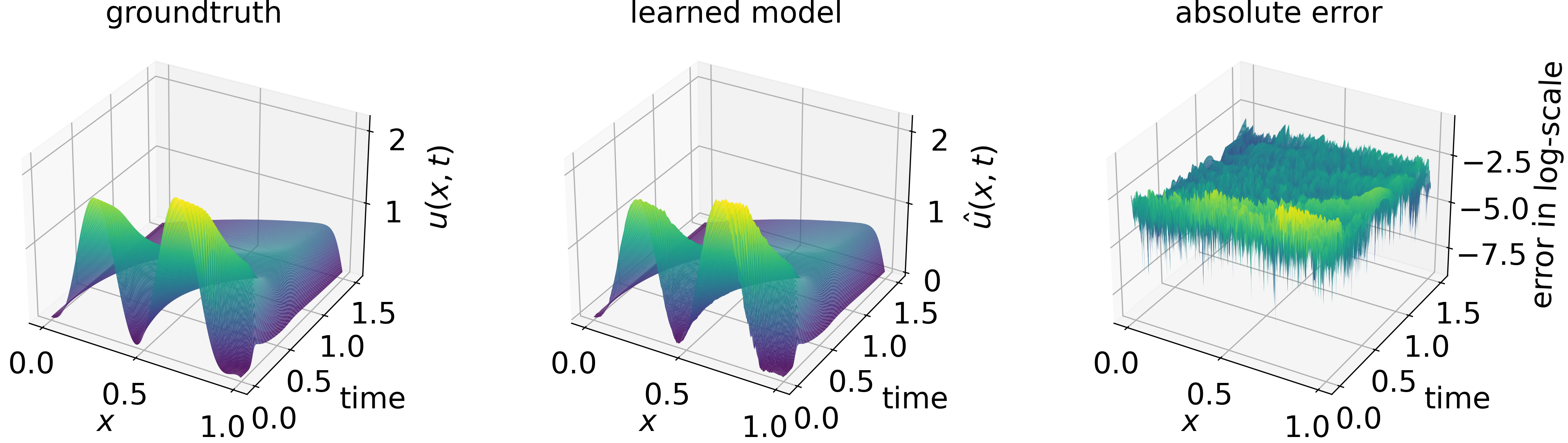}
	\includegraphics[width = 0.75\linewidth,angle = 0]{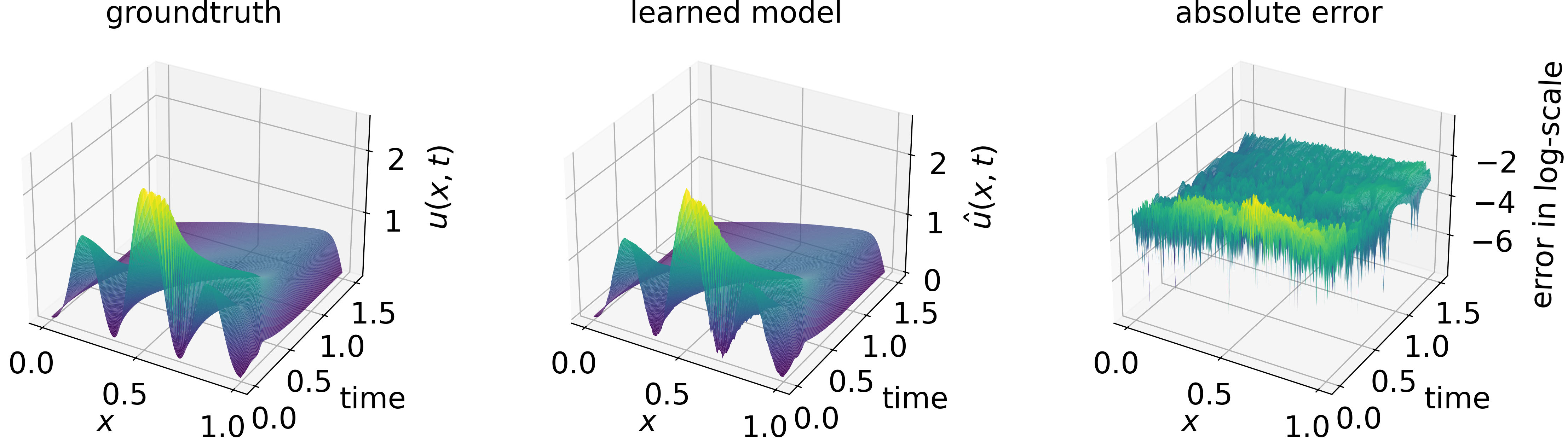}
	\caption{Burgers' equation: A comparison of the solutions using  \texttt{quad-embeds-AE} with the ground truth for two test cases.}
	\label{fig:burgers_time_domain}
\end{figure}

Unlike the previous two examples, this example is high-dimensional. However, it is well-known that high-dimensional systems often evolve in a low-dimensional subspace. Using this hypothesis, one can project a high-dimensional nonlinear system into a low-dimensional subspace to obtain a low-dimensional nonlinear system. Additionally, we know that smooth nonlinear systems can be written as quadratic systems. Combining these two philosophies for high-dimensional nonlinear systems, we aim to learn a low-dimensional representation so that a quadratic system can govern its dynamics. To determine a low-dimensional representation, we make use of an autoencoder consisting of convolutional layers. Its detailed design is discussed in \ref{sec:appendix}. We fix the latent representation dimension to four and train the parameters of autoencoders and quadratic systems to describe the dynamics of the latent representation. We denote this approach as \texttt{quad-embeds-conv}. The performance of the proposed method for two test cases out of four is illustrated in \Cref{fig:burgers_time_domain}, where we notice a faithful recovery of the dynamics. However, we notice small perturbations in the learned solutions, which, we believe, can be reduced with a more powerful decoder design, potentially using residual connections, but this would come with more computational expenses and data requirements. 

\begin{figure}[tb]
	\centering
	\includegraphics[width = 0.75\linewidth,angle = 0]{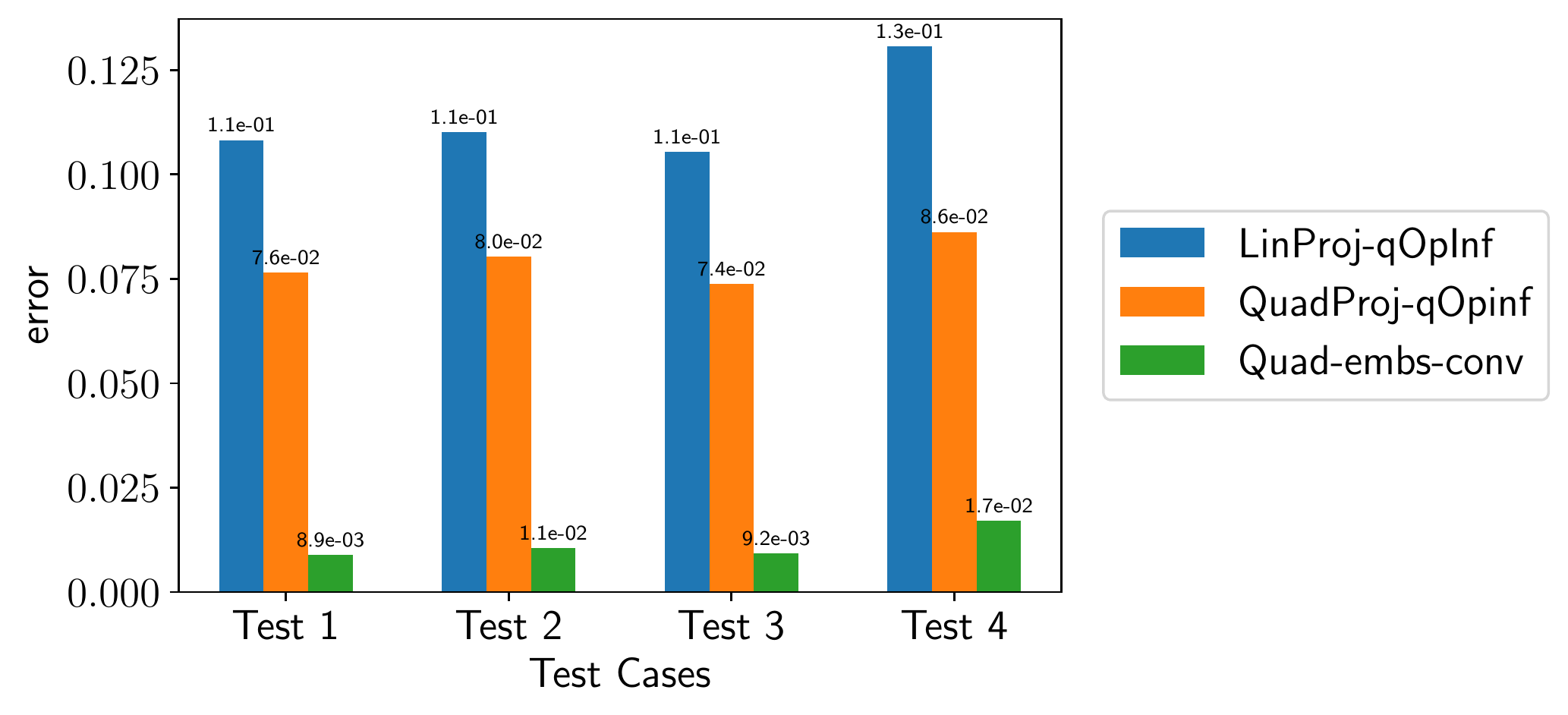}
	\caption{Burgers' equation: A comparison of the solutions using  \texttt{quad-embeds-AE} with the ground truth for two test cases.}
	\label{fig:burgers_comparison}
\end{figure}

We compare the performance of \texttt{quad-embeds-conv} with operator inference \cite{morPehW16}. For this, we first project the high-dimensional data using the four most dominant POD modes and learn a quadratic model. The solution of the full spatial grid is obtained by re-projecting the low-dimensional solutions using the same POD basis. We denote this approach by \texttt{LinProj-qOpInf}. The second method for a comparison is based on the work \cite{geelen2023operator}, which is inspired by the quadratic manifold work \cite{barnett2022quadratic}. The method in \cite{geelen2023operator},  in principle, aims to learn a quadratic model using  POD coordinates that are obtained by projecting high-dimensional training data using the dominant POD modes, and an approximate solution on the full spatial solution is obtained by a quadratic ansatz. We refer to this approach as \texttt{QuadProj-qOpInf}. For  \texttt{LinProj-qOpInf} and \texttt{QuadProj-qOpInf}, we consider the four most dominant POD modes. Next, we present a comparison of these three methods for four test cases using the following measure:
\begin{equation}\label{eq:error_measure_relative}
	\cE(f) =   \dfrac{\left\|\bX_{\texttt{ground-truth}}^{(f)} - \bX_{\texttt{learned}}^{(f)}  \right\|_2}{\left\|\bX_{\texttt{ground-truth}}^{(f)}   \right\|_2},
\end{equation}
where $\bX_{\texttt{ground-truth}}^{(f)}$ and  $\bX_{\texttt{learned}}^{(f)}$, respectively, are solutions using the ground truth model and the learned models for a given frequency $f$ that defines a test initial condition. The result is shown in \Cref{fig:burgers_comparison}. We notice that \texttt{QuadProj-qOpInf} performs better as compared to \texttt{LinProj-qOpInf} as expected. But \texttt{quad-embeds-conv} out-performances both approaches by a margin. 

\section{Discussions}
In this work, we have discussed a unified representation for nonlinear systems, the so-called quadratic embeddings. This idea stems from the fact that smooth nonlinear dynamical systems can be written as quadratic systems in appropriate lifted coordinates. While it is possible to manually design lifted coordinates for specific analytical expressions of nonlinear systems, this becomes challenging in a data-driven context where the goal is to learn dynamical models directly from data.
To address this challenge, in this work, we have utilized the powerful approximation capabilities of deep neural networks to learn lifted coordinates so that a quadratic model can describe its dynamics. We have proposed the usage of an autoencoder design that aims to learn lifted coordinates so that a quadratic system can describe its dynamics. As a result, we obtain parsimonious models since the learned differential equations have a simple quadratic analytical expression, which can be exploited in engineering design. We have also discussed asymptotic stability for the dynamics of the lifted coordinates via parameterization of globally stable quadratic systems. We have also discussed an extension to high-dimensional nonlinear systems. We have compared our approach with the universal linear embedding method and operator inference, demonstrating the efficiency of our proposed method.

This work opens several promising avenues for further research. Since the autoencoder that determines lifted variables is associated with deep neural networks, it demands a diverse dataset and intense computational resources. Although recent advancements in GPUs and efficient deep learning libraries such as TensorFlow \cite{tensorflow2015-whitepaper} and PyTorch \cite{NEURIPS2019_9015} can allow overcoming computational difficulties, the interpretability and generalizability of these networks still need further investigation. Often, deep learning models are good interpolating models, but their extrapolation capabilities are questionable. Therefore, it is essential to incorporate physical knowledge and any available prior knowledge about the origin of the data in the learning process. As a result, this is not only expected to improve the interpretability and generalizability of the learned model but also may extrapolate better outside the training regime. It may also reduce the data required for training as scarcity of it can be compensated by embedding physical laws. Moreover, one of the critical hyper-parameters in our approach is the dimension of the lifted coordinate system. Finding a minimal dimension of it is desirable while keeping the prescribed accuracy. It would be worthwhile to develop a suitable automatic approach to determine the dimension of the lifted coordinate system. Additionally, using physical and domain knowledge in deep learning frameworks can enhance interpretability and provide some physics-informed quadratic embeddings, i.e., quadratic embeddings for Hamiltonian systems. The proposed framework can be extended to other classes, such as systems with parameters and control. Special treatment is required for noisy measurements, for which ideas proposed in \cite{rudy2019deep,morGoyB21b,morGoyB22} can be combined with our methodology. In the future, we will apply our approach to more challenging and important applications in science and technology to construct models where first-principle modeling remains challenging, and the learned models can enhance the engineering design process.

\addcontentsline{toc}{section}{References}
\bibliographystyle{ieeetr}
\bibliography{mybib}

\begin{thebibliography}{10}

\bibitem{crutchfield1987equations}
J.~P. Crutchfield and B.~S. McNamara, ``Equations of motion from a data
  series,'' {\em Complex Sys.}, vol.~1, no.~121, pp.~417--452, 1987.

\bibitem{bongard2007automated}
J.~Bongard and H.~Lipson, ``Automated reverse engineering of nonlinear
  dynamical systems,'' {\em Proc. Nat. Acad. Sci. U.S.A.}, vol.~104, no.~24,
  pp.~9943--9948, 2007.

\bibitem{yao2007modeling}
C.~Yao and E.~M. Bollt, ``Modeling and nonlinear parameter estimation with
  {K}ronecker product representation for coupled oscillators and spatiotemporal
  systems,'' {\em Physica D: Nonlinear Phenomena}, vol.~227, no.~1, pp.~78--99,
  2007.

\bibitem{schmidt2011automated}
M.~D. Schmidt, R.~R. Vallabhajosyula, J.~W. Jenkins, J.~E. Hood, A.~S. Soni,
  J.~P. Wikswo, and H.~Lipson, ``Automated refinement and inference of
  analytical models for metabolic networks,'' {\em Phy. Biology}, vol.~8,
  no.~5, p.~055011, 2011.

\bibitem{rowley2009spectral}
C.~W. Rowley, I.~Mezi\'c, S.~Bagheri, P.~Schlatter, and D.~Henningson,
  ``Spectral analysis of nonlinear flows,'' {\em J. Fluid Mech.}, vol.~641,
  no.~1, pp.~115--127, 2009.

\bibitem{schmid2010dynamic}
P.~J. Schmid, ``Dynamic mode decomposition of numerical and experimental
  data,'' {\em J. Fluid Mech.}, vol.~656, pp.~5--28, 2010.

\bibitem{morBenGW15}
P.~Benner, S.~Gugercin, and K.~Willcox, ``A survey of projection-based model
  reduction methods for parametric dynamical systems,'' {\em {SIAM} Rev.},
  vol.~57, no.~4, pp.~483--531, 2015.

\bibitem{morPehW16}
B.~Peherstorfer and K.~Willcox, ``Data-driven operator inference for
  nonintrusive projection-based model reduction,'' {\em Comp. Meth. Appl. Mech.
  Eng.}, vol.~306, pp.~196--215, 2016.

\bibitem{rudy2017data}
S.~H. Rudy, S.~L. Brunton, J.~L. Proctor, and J.~N. Kutz, ``Data-driven
  discovery of partial differential equations,'' {\em Sci. Adv.}, vol.~3,
  no.~4, p.~e1602614, 2017.

\bibitem{takeishi2017learning}
N.~Takeishi, Y.~Kawahara, and T.~Yairi, ``Learning {K}oopman invariant
  subspaces for dynamic mode decomposition,'' {\em Adv. Neural Inform.
  Processing Systems}, vol.~30, pp.~1130--1140, 2017.

\bibitem{morGoyB22a}
P.~Goyal and P.~Benner, ``Discovery of nonlinear dynamical systems using a
  {R}unge-{K}utta inspired dictionary-based sparse regression approach,'' {\em
  Proc. Royal Society A: Mathematical, Physical and Engineering Sciences},
  vol.~478, no.~2262, p.~20210883, 2022.

\bibitem{morMayA07}
A.~J. Mayo and A.~C. Antoulas, ``A framework for the solution of the
  generalized realization problem,'' {\em Linear Algebra Appl.}, vol.~425,
  no.~2--3, pp.~634--662, 2007.

\bibitem{morDrmGB15a}
Z.~Drma{\v{c}}, S.~Gugercin, and C.~Beattie, ``Vector fitting for matrix-valued
  rational approximation,'' {\em {SIAM} J. Sci. Comput.}, vol.~37, no.~5,
  pp.~A2346--A2379, 2015.

\bibitem{nakatsukasa2018aaa}
Y.~Nakatsukasa, O.~S{\`e}te, and L.~N. Trefethen, ``The {AAA} algorithm for
  rational approximation,'' {\em {SIAM} J. Sci. Comput.}, vol.~40, no.~3,
  pp.~A1494--A1522, 2018.

\bibitem{rico1994continuous}
R.~Rico-Martinez, J.~Anderson, and I.~Kevrekidis, ``Continuous-time nonlinear
  signal processing: a neural network based approach for gray box
  identification,'' in {\em Proc. IEEE Workshop on Neural Networks for Signal
  Processing}, pp.~596--605, IEEE, 1994.

\bibitem{rico1995nonlinear}
R.~Rico-Martinez, I.~Kevrekidis, and K.~Krischer, ``Nonlinear system
  identification using neural networks: dynamics and instabilities,'' {\em
  Neural Networks for Chemical Engineers}, pp.~409--442, 1995.

\bibitem{gonzalez1998identification}
R.~Gonzalez-Garcia, R.~Rico-Martinez, and I.~Kevrekidis, ``Identification of
  distributed parameter systems: {A} neural net based approach,'' {\em
  Computers \& Chemical Engrg.}, vol.~22, pp.~S965--S968, 1998.

\bibitem{mardt2018vampnets}
A.~Mardt, L.~Pasquali, H.~Wu, and F.~No{\'e}, ``{VAMPnets} for deep learning of
  molecular kinetics,'' {\em Nature Commun.}, vol.~9, no.~1, pp.~1--11, 2018.

\bibitem{vlachas2018data}
P.~R. Vlachas, W.~Byeon, Z.~Y. Wan, T.~P. Sapsis, and P.~Koumoutsakos,
  ``Data-driven forecasting of high-dimensional chaotic systems with long
  short-term memory networks,'' {\em Proc. Royal Society A: Mathematical,
  Physical and Engineering Sciences}, vol.~474, no.~2213, p.~20170844, 2018.

\bibitem{champion2019data}
K.~Champion, B.~Lusch, J.~N. Kutz, and S.~L. Brunton, ``Data-driven discovery
  of coordinates and governing equations,'' {\em Proc. Nat. Acad. Sci. U.S.A.},
  vol.~116, no.~45, pp.~22445--22451, 2019.

\bibitem{lusch2018deep}
B.~Lusch, J.~N. Kutz, and S.~L. Brunton, ``Deep learning for universal linear
  embeddings of nonlinear dynamics,'' {\em Nature Commun.}, vol.~9, no.~1,
  pp.~1--10, 2018.

\bibitem{chen2018neural}
R.~T. Chen, Y.~Rubanova, J.~Bettencourt, and D.~K. Duvenaud, ``Neural ordinary
  differential equations,'' in {\em Advances Neural Inform. Processing Sys.},
  pp.~6571--6583, 2018.

\bibitem{morGoyB21b}
P.~Goyal and P.~Benner, ``Learning dynamics from noisy measurements using deep
  learning with a {R}unge-{K}utta constraint,'' {\em Workshop paper at the
  Symbiosis of Deep Learning and Differential Equations -- NeurIPS.~}, 2021.
\newblock Avaiable at https://openreview.net/forum?id=G5i2aj7v7i.

\bibitem{ogata2010modern}
K.~Ogata {\em et~al.}, {\em Modern Control Engineering}, vol.~5.
\newblock Prentice Hall Upper Saddle River, NJ, 2010.

\bibitem{aastrom2021feedback}
K.~J. {\AA}str{\"o}m and R.~M. Murray, {\em Feedback Systems: An Introduction
  for Scientists and Engineers}.
\newblock Princeton University Press, 2021.

\bibitem{brunton2022data}
S.~L. Brunton and J.~N. Kutz, {\em Data-Driven Science and Engineering: Machine
  Learning, Dynamical Systems, and Control}.
\newblock Cambridge University Press, 2022.

\bibitem{koopman1931hamiltonian}
B.~O. Koopman, ``Hamiltonian systems and transformation in {H}ilbert space,''
  {\em Proc. Nat. Acad. Sci. U.S.A.}, vol.~17, no.~5, p.~315, 1931.

\bibitem{williams2015data}
M.~O. Williams, I.~G. Kevrekidis, and C.~W. Rowley, ``A data--driven
  approximation of the {K}oopman operator: Extending dynamic mode
  decomposition,'' {\em J. Nonlinear Sci.}, vol.~25, no.~6, pp.~1307--1346,
  2015.

\bibitem{morKutBBetal16}
J.~N. Kutz, S.~L. Brunton, B.~W. Brunton, and J.~L. Proctor, {\em Dynamic Mode
  Decomposition: Data-Driven Modeling of Complex Systems}.
\newblock Philadelphia, USA: Society of Industrial and Applied Mathematics,
  2016.

\bibitem{li2017extended}
Q.~Li, F.~Dietrich, E.~M. Bollt, and I.~G. Kevrekidis, ``Extended dynamic mode
  decomposition with dictionary learning: A data-driven adaptive spectral
  decomposition of the {K}oopman operator,'' {\em Chaos: An Interdisciplinary
  Journal of Nonlinear Science}, vol.~27, no.~10, pp.~103--111, 2017.

\bibitem{morBenHM18}
P.~Benner, C.~Himpe, and T.~Mitchell, ``On reduced input-output dynamic mode
  decomposition,'' {\em Adv. Comput. Math.}, vol.~44, no.~6, pp.~1821--1844,
  2018.

\bibitem{mezic2020spectrum}
I.~Mezi{\'c}, ``Spectrum of the {K}oopman operator, spectral expansions in
  functional spaces, and state-space geometry,'' {\em J. Nonlin. Sci.},
  vol.~30, no.~5, pp.~2091--2145, 2020.

\bibitem{savageau1987recasting}
M.~A. Savageau and E.~O. Voit, ``Recasting nonlinear differential equations as
  {S}-systems: {A} canonical nonlinear form,'' {\em Mathematical Biosciences},
  vol.~87, no.~1, pp.~83--115, 1987.

\bibitem{papachristodoulou2005analysis}
A.~Papachristodoulou and S.~Prajna, ``Analysis of non-polynomial systems using
  the sum of squares decomposition,'' in {\em Positive Polynomials in Control}
  (A.~Henrion, Didierand~Garulli, ed.), pp.~23--43, Springer, 2005.

\bibitem{morGu11}
C.~Gu, ``{QLMOR}: A projection-based nonlinear model order reduction approach
  using quadratic-linear representation of nonlinear systems,'' {\em IEEE
  Trans. Comput. Aided Des. Integr. Circuits. Syst.}, vol.~30, no.~9,
  pp.~1307--1320, 2011.

\bibitem{morBenB15}
P.~Benner and T.~Breiten, ``Two-sided projection methods for nonlinear model
  order reduction,'' {\em {SIAM} J. Sci. Comput.}, vol.~37, no.~2,
  pp.~B239--B260, 2015.

\bibitem{qian2020lift}
E.~Qian, B.~Kramer, B.~Peherstorfer, and K.~Willcox, ``Lift \& learn:
  Physics-informed machine learning for large-scale nonlinear dynamical
  systems,'' {\em Physica D: Nonlinear Phenomena}, vol.~406, p.~132401, 2020.

\bibitem{morSchVR08}
W.~H.~A. Schilders, H.~A. {van der Vorst}, and J.~Rommes, {\em Model Order
  Reduction: Theory, Research Aspects and Applications}.
\newblock Berlin, Heidelberg: Springer-Verlag, 2008.

\bibitem{lee2020model}
K.~Lee and K.~T. Carlberg, ``Model reduction of dynamical systems on nonlinear
  manifolds using deep convolutional autoencoders,'' {\em J. Comput. Phys.},
  vol.~404, p.~108973, 2020.

\bibitem{barnett2022quadratic}
J.~Barnett and C.~Farhat, ``Quadratic approximation manifold for mitigating the
  {K}olmogorov barrier in nonlinear projection-based model order reduction,''
  {\em J. Comput. Phys.}, vol.~464, p.~111348, 2022.

\bibitem{geelen2023operator}
R.~Geelen, S.~Wright, and K.~Willcox, ``Operator inference for non-intrusive
  model reduction with quadratic manifolds,'' {\em Comp. Meth. Appl. Mech.
  Eng.}, vol.~403, p.~115717, 2023.

\bibitem{mccormick1976computability}
G.~P. McCormick, ``Computability of global solutions to factorable nonconvex
  programs: {Part I—Convex} underestimating problems,'' {\em Mathematical
  Programming}, vol.~10, no.~1, pp.~147--175, 1976.

\bibitem{amato2007state}
F.~Amato, R.~Ambrosino, M.~Ariola, C.~Cosentino, and A.~Merola, ``State
  feedback control of nonlinear quadratic systems,'' in {\em 46th IEEE Conf.
  Decision and Control}, pp.~1699--1703, IEEE, 2007.

\bibitem{tognetti2021output}
E.~S. Tognetti, M.~Jungers, and T.~R. Calliero, ``Output feedback control for
  quadratic systems: A {L}yapunov function approach,'' {\em Internat. J. Robust
  Nonlinear Control}, vol.~31, no.~17, pp.~8373--8389, 2021.

\bibitem{goyal2023guaranteed}
P.~Goyal, I.~P. Duff, and P.~Benner, ``Guaranteed stable quadratic models and
  their applications in {SINDy} and operator inference,'' {\em arXiv preprint
  arXiv:2308.13819}, 2023.

\bibitem{nutku1990hamiltonian}
Y.~Nutku, ``Hamiltonian structure of the {L}otka-{V}olterra equations,'' {\em
  Physics Letters A}, vol.~145, no.~1, pp.~27--28, 1990.

\bibitem{tensorflow2015-whitepaper}
M.~Abadi {\em et~al.}, ``{TensorFlow}: Large-scale machine learning on
  heterogeneous systems,'' 2015.
\newblock Software available from tensorflow.org.

\bibitem{NEURIPS2019_9015}
A.~Paszke {\em et~al.}, ``{PyTorch}: An imperative style, high-performance deep
  learning library,'' in {\em Advances in Neural Inform. Processing Systems}
  (H.~Wallach, H.~Larochelle, A.~Beygelzimer, F.~d\textquotesingle
  Alch\'{e}-Buc, E.~Fox, and R.~Garnett, eds.), vol.~32, pp.~8024--8035, 2019.

\bibitem{rudy2019deep}
S.~H. Rudy, J.~N. Kutz, and S.~L. Brunton, ``Deep learning of dynamics and
  signal-noise decomposition with time-stepping constraints,'' {\em J. Comput.
  Phys.}, vol.~396, pp.~483--506, 2019.

\bibitem{morGoyB22}
P.~Goyal and P.~Benner, ``Neural ordinary differential equations with irregular
  and noisy data,'' {\em Royal Society Open Sci.}, vol.~10, no.~7, p.~221475,
  2023.

\end{thebibliography}

\appendix
\section{Implementation details}\label{sec:appendix}
\begin{figure}[tb]
	\centering
	\includegraphics[width=1\linewidth,height=0.15\textheight]{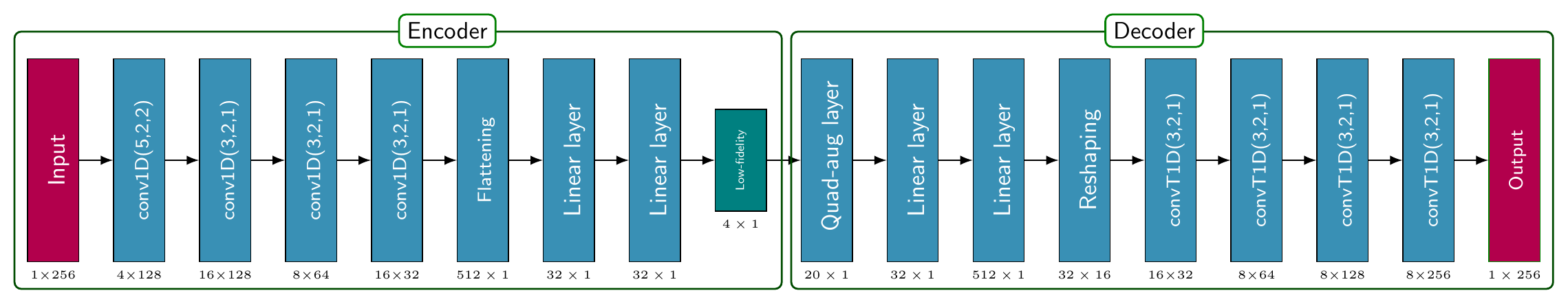}
	\caption{The figure illustrates the design of the autoencoder for the Burgers' example. The autoencoder contains the encoder and decoder parts. In the figure, conv1D$(k,s,p)$ denotes a 1D convolution layer with kernel size $k$, stride size $s$, padding size $p$, and similarly, convT1D$(k,s,p)$ is a 1D transpose convolution layer with transpose kernel size $k$, stride size $s$, padding size $p$. We show the size of the output block below each layer. Furthermore, the decoder part has a customized layer, namely~\texttt{Quad-aug layer}; it augments input and Kronecker product of inputs, which is then passed to the later layers.}
	\label{fig:dcaschematic}
\end{figure}

Here, we provide the necessary details used in our experiments. All the experiments are carried out a machine with an \texttt{Intel\textsuperscript{\tiny\textcopyright} Core\textsuperscript{\tiny TM} i5-12600K} CPU and \texttt{NVIDIA RTX\textsuperscript{\tiny TM} A4000(32GB)} GPU. We have generated training and testing data for random initial conditions using the function \texttt{solve\_ivp} with default parameter settings from \texttt{scipy} library. 
For nonlinear pendulum and dissipative Lokta-Volterra examples, we have used an autoencoder design based on multi-layer fully connected neural networks with skip connections, and for Burgers' example, the used autoencoder design is shown in \Cref{fig:dcaschematic}.
\Cref{tab:hyperparameters_lowdimensional} contains all other necessary hyper-parameters for the illustrative examples. 
Based on a few trials, we set $(\lambda_1,\lambda_2,\lambda_3)$ for each example as shown in \Cref{tab:hyperparameters_lowdimensional}. Note, we set $\lambda_2$ to zero to avoid additional computational of derivatives through decoder. However, we believe that determining a good balance for these different losses using a proper cross-validation can improve the performance of our approach.
For training, we have also utilized a decaying learning rate, for which we have reduced our learning rate by $\tfrac{1}{10}$ after every $M$ number of epoch (the value of $M$ is given in \Cref{tab:hyperparameters_lowdimensional}). Furthermore, we have utilized a hard-pruning feature in our training. We initialize parameters of autoencoder using by default PyTorch (version 1.13) settings, and the matrices $\{\bA,\bH\}$ are initialized as suggested in \cite{goyal2023guaranteed}.

\begin{table}[tb]
	\renewcommand{\arraystretch}{1.25}
	\caption{The table contains all the hyper-parameters to learn the autoencoder parameters and matrices defining lifted-coordinate dynamics. }
	\label{tab:hyperparameters_lowdimensional}
	\begin{tabular}{|c|c|c|c|}
		\hline
		Parameters                                                              & \begin{tabular}[c]{@{}c@{}}Pendulum \\ example\end{tabular} & \begin{tabular}[c]{@{}c@{}}Dissipative Lotka-Volterra \\ example\end{tabular} & \begin{tabular}[c]{@{}c@{}}Burgers' \\ example\end{tabular} \\ \hline
		Number of neurons  & $8$   & $16$       & N/A \\ \hline
		\begin{tabular}[c]{@{}c@{}}Lifted coordinate \\ system dimension\end{tabular} & 3                                                           & 3 & 4                                                                                                                                   \\ \hline
		Learning rate   & $3\cdot 10^{-3}$    &$3\cdot 10^{-3}$ &$5\cdot 10^{-3}$  \\ \hline
		Batch size                                                             & $32$                                                    & $64$    & $64$                                                                         \\ \hline
		Activation function                                                    & \texttt{silu}                              & \texttt{silu}    & N/A                                                                  \\ \hline
		Weight decay                                                           & $10^{-5}$                                                   &     $10^{-5}$   &     $10^{-5}$                                                \\ \hline
		Epochs                                                                 & 4000                                                        &    4000    & 400                                                   \\ \hline
		$M$                                                                 & 1500                                                        &    1500    & 150                                                   \\ \hline
		$(\lambda_1, \lambda_2,\lambda_3)$                                                                 & (1, 0, 1)                                                        &    (1,0,1)      &  (10,0,1)                                               \\ \hline
	\end{tabular}
\end{table}

\end{document}